\documentclass[11pt]{article}
\usepackage{latexsym}
\usepackage{amsmath}
\usepackage{cite}
\usepackage{amsthm,color,bm}
\usepackage{amssymb}
\usepackage{mathrsfs}
\usepackage{indentfirst}
\usepackage{times,cases}
\usepackage{dsfont}

\usepackage{lscape}
\newtheorem{theorem}{\noindent Theorem}[section]
\newtheorem{example}{\noindent Example}[section]
\newtheorem{proposition}[theorem]{\noindent Proposition}
\newtheorem{definition}[theorem]{\noindent Definition}
\newtheorem{lemma}[theorem]{\noindent Lemma}
\newtheorem{remark}[theorem]{\noindent Remark}
\newtheorem{corollary}[theorem]{\noindent Corollary}
\numberwithin{figure}{section}
\numberwithin{equation}{section}

\renewcommand{\theequation}{\thesection.\arabic{equation}}

\newcommand{\fR}{\mathfrak{R}}

\def\Z{\mathbb{Z}}
\def\R{\mathbb{R}}

\def\C{\mathbb{C}}

\def\N{\mathbb{N}}
\def\1{\mathds{1}}

\def\disp{\displaystyle}

\def\bc{\begin{center}}
\def\ec{\end{center}}
\def\be{\begin{equation}}
\def\ee{\end{equation}}
\def\bea{\begin{eqnarray}}
\def\eea{\end{eqnarray}}
\def\ba{\begin{array}}
\def\ea{\end{array}}
\def\benu{\begin{enumerate}}
\def\eenu{\end{enumerate}}
\def\bt{\begin{theorem}}
\def\et{\end{theorem}}
\def\bl{\begin{lemma}}
\def\el{\end{lemma}}
\def\bco{\begin{corollary}}
\def\eco{\end{corollary}}
\def\bn{\begin{numcases}}
\def\en{\end{numcases}}
\def\br{\begin{remark}}
\def\er{\end{remark}}
\def\bd{\begin{definition}}
\def\ed{\end{definition}}
\def\bp{\begin{proposition}}
\def\ep{\end{proposition}}
\def\bo{\begin{proof}}
\def\eo{\end{proof}}
\def\bx{\begin{example}}
\def\ex{\end{example}}
\def\bal{\begin{align}}
\def\eal{\end{align}}

\def\al{\alpha}

\def\lam{\lambda} 

\def\ve{\varepsilon}

\def\vp{\varphi}

\def\~{\widetilde}

\def\ra{\rightarrow}
\def\Ra{\Rightarrow}
\def\Lra{\Leftrightarrow}

\def\8{\infty}

\def\mb{\mbox}

\def\sm{\setminus}

\def\hs{\hspace{0.4cm}}
\def\Vs{\vskip10pt}
\def\vs{\vskip5pt}

\def\({\left(}
\def\){\right)}

\begin{document}


\begin{center}
    {\large \bf New high-dimensional generalizations of Nesbitt's
    inequality and relative applications}
\vspace{0.5cm}\\
{Junfeng Zhang,\quad Jintao Wang$^*$}\\\vspace{0.3cm}

{\small Department of Mathematics, Wenzhou University, Wenzhou 325035, China}
\end{center}


\renewcommand{\theequation}{\arabic{section}.\arabic{equation}}
\numberwithin{equation}{section}


\begin{abstract}
Two kinds of novel generalizations of Nesbitt's inequality are explored in various cases
regarding dimensions and parameters in this article.
Some other cases are also discussed elaborately by using the semiconcave-semiconvex theorem.
The general inequalities are then employed to deduce some alternate inequalities and
mathematical competition questions.
At last, a relation about Hurwitz-Lerch zeta functions is obtained.
\vs

\noindent\textbf{Keywords:} Nesbitt's inequality, Jensen's inequality, Chepyshev's inequality, semiconcave-semiconvex theorems, Hurwitz-Lerch zeta functions.
\vs

\noindent{\bf AMS Subject Classification 2010:}\, 26D15, 26D10, 11M35

\end{abstract}

\vspace{-1 cm}

\footnote[0]{\hspace*{-7.4mm}
$^{*}$ Corresponding author.\\
E-mail address: wangjt@wzu.edu.cn (J.T. Wang),\hs 21211357329@stu.wzu.edu.cn (J. F. Zhang).}

\section{Introduction}
Inequalities play a significant and fundamental role in the development of modern science,
technology and education (\cite{K04}).
As an ancient Chinese proverb goes, ``A very tiny difference within a millimeter can lead to
an error of more than thousands of miles", which is just like a fatal tornado caused by a butterfly's
flapping wings.
Since it is impossible to measure and constrain the real things in the absolute sense,
the most important issue we have to face is how to estimate and ascertain
the terrible unknown outcomes.
In this process, inequalities have showcased extraordinary application value
(\cite{FY22,LW21-1,LW25,LLJ21,WJ24,WL21,WLD21,WLYJ21,WPL25,WWL22,WZL23,WZL24,WQB13}).

In area of education, inequalities are of particular effectiveness to practice and test
the intelligence of students in high school (\cite{C13,HL16,J12,JC21,JG23,L06,W09,W14,ZW11}).
Typical ones of such inequalities include alternate inequalities, mean value inequalities,
and Radon's inequality.
Amongst these inequalities, Nesbitt's inequality (see \cite{N03})
$$\frac{x}{y+z}+\frac{y}{z+x}+\frac{z}{x+y}\geqslant\frac32$$
has been known as a famous one and generalized to different forms since 1903.
In recent two decades, increasing attention has been paid to generalizations and
relative applications of Nesbitt's inequality.
Bencze et al in \cite{BZ08,BP11} gave one kind of generalization with weights
and refinements of Nesbitt's inequality.
Batinetu-Giurgiu and Stanciu presented some concrete examples of generalizations with weights
and analogous form of Nesbitt's inequality in \cite{BS12,BS13-1,BS13-2}.
An iconic generalization of Nesbitt's inequality was a high-dimensional version given
by Wang in \cite{WHJ15} and read as follows:
\be\label{1.1}\sum_{i=1}^n\frac{a_i^m}{ts-a_i^p}\geqslant\frac{\disp\sum_{i=1}^na_i^{m-p}}{nt-1},\ee
where $a_1,\cdots,a_n>0$, $m,n,p,t\in\N^+$, $m\geqslant p$, $n\geqslant2$
and $\disp s=\sum_{i=1}^na_i^p$.
Chu, Jiang et al also generalized Nesbitt's inequality on dimensions and (integer) powers in
\cite{C17,JC21,JG23}.
The Nesbitt's inequality was also concerned with in the study of other inequalities
(\cite{BS13-2,J12,M11,WW09}).
What is even more interesting is that Nesbitt's inequality can be also applied
to other fields such as theories of matrices and numbers (\cite{CS22,WQB13}).

In this article, we further develop and generalize Nesbitt's inequality with more parameters
and high dimensions in different forms.
The newly generalized versions of Nesbitt's inequality cover most generalized versions given before
and even include the situations that derive inverse inequalities.
Specifically, we consider the algebraic expression
\be\label{1.2}\sum_{i=1}^n\frac{a_i^m}{(ts-ra_i^p)^\beta},\hs\mb{where}\hs s=\sum_{i=1}^na_i^p,\ee
and compare \eqref{1.2} with
\be\label{1.3}\frac{n^{\beta+1-\frac mp}}{(nt-r)^\beta}\(\sum_{i=1}^na_i^p\)^{\frac mp-\beta}\hs\mb{and}
\hs\frac{1}{(nt-r)^{\beta}}\sum_{i=1}^na_i^{m-\beta p},\ee
where $n\in\N^+$, $a_1,\cdots,a_n>0$, $m,p,\beta,t,r\in\R$ with $t\geqslant0$ and $ts>ra_i^p$
for all $i=1,\cdots,n$.
The inequality \eqref{1.1} is a simple relation of \eqref{1.2} and (the second expression of)
\eqref{1.3} for the case when $\beta=1$.

Our main goal in this article is to study the relation between \eqref{1.2} and \eqref{1.3},
which differs greatly in different cases.
To compare \eqref{1.2} and the first algebraic expression of \eqref{1.3}
in Theorem \ref{th3.1}, the Jensen's inequality is a powerful tool,
and a generalized version (Theorem \ref{GWSIP}) of Radon's inequality is also of great help.
To determine the relation of \eqref{1.2} and the second one of \eqref{1.3},
we employ Theorem \ref{th3.1}, rearragement inequality, Chepyshev's inequality and
Jensen's inequlity and give definitive results in different
cases in Theorem \ref{th3.2} and \ref{th3.3}.
The inequality consequences proved above do not cover all the cases.
For other cases that guarantee the inequlities, a useful theorem --- Semiconcave-semiconvex theorem
from \cite{H11} is rather effective.

The newly generalized Nesbitt's inequalities can be applied to
prove many alternate inequalities in different forms regarding dimensions, parameters and exponents.
In particular, some competitive contest questions, including international mathematical Olympiad
(IMO for short) questions, can be easily obtained only by picking certain parameters in the generalized inequalities.

At last, we also consider the applications of the obtained inequalities in the study of Hurwitz-Lerch functions.
In \cite{WQB13}, Wang obtained the minimum value related to Riemann's and Hurwitz's
zeta function by using his main inequality
$$\frac{\disp\(\sum_{k=1}^np_kx_k\)^{\al}}{\disp\(M-\sum_{k=1}^np_kx_k\)^\beta}\leqslant
\sum_{k=1}^n\frac{p_kx_k^\al}{(M-x_k)^\beta},$$
where $\al\geqslant1$, $\beta\geqslant0$, $0<x_k<M<+\8$, $p_k\in[0,1]$, $k=1,\cdots,n$
with $p_1+\cdots+p_n=1$.
In our work, we further study the relation of different Hurwitz-Lerch functions by using our generalized inequalities.
We not only generalize the result of \cite{WQB13}, but also obtain a new inverse relation.

The remainder of this article is organized as follows.
In Section 2, some necessary inequalities are presented for the following argument.
In Section 3, the main theorems are proved and some examples of other cases are given for clarity.
In Section 4, we apply the main theorems to some inequality problems and competition questions.
In Section 5, we apply the main theorems to obtaining some relations about
different Hurwitz-Lerch functions.

\section{Preliminaries}\label{s2}
In this section, we present some necessary basic inequalities.

First we recall the \textbf{Rearrangement Inequality}.
Let $a_i$, $b_i\in\R$ ($1\leqslant i\leqslant n$) with
\be\label{2.1}a_1\leqslant a_2\leqslant\cdots\leqslant a_n\hs\mb{and}\hs
b_1\leqslant b_2\leqslant\cdots\leqslant b_n,\ee
and $\{c_i\}_{1\leqslant i\leqslant n}$ be
a rearrangement of $\{b_i\}_{1\leqslant i\leqslant n}$.
Then it holds that
$$\sum_{i=1}^na_ib_{n+1-i}\leqslant\sum_{i=1}^na_ic_i
\leqslant\sum_{i=1}^na_ib_i.$$
Applying the rearrangement inequality stated above,
one can easily obtain the \textbf{Chepyshev's inequality}:
for $\{a_i\}_{1\leqslant i\leqslant n}$,
$\{b_i\}_{1\leqslant i\leqslant n}$ given in \eqref{2.1},
it holds that
$$\sum_{i=1}^na_ib_i\geqslant\frac{1}{n}\sum_{i=1}^na_i\cdot\sum_{i=1}^nb_i
\geqslant\sum_{i=1}^na_ib_{n+1-i}.$$

We then recall the famous \textbf{Jensen's inequality}.
Let $I\subset\R$ be an interval, $\vp:I\ra\R$ a convex function,
$\psi:I\ra\R$ a concave one,
then for each $n\in\N$, $x_1$, $\cdots$, $x_n\in I$ and positive weights $\lam_1$,
$\cdots$, $\lam_n$ with $\lam_1+\cdots+\lam_n=1$, the following inequalities hold:
$$\vp\(\sum_{i=1}^n\lam_ix_i\)\leqslant\sum_{i=1}^n\lam_i\vp(x_i)\hs\mb{and}\hs
\psi\(\sum_{i=1}^n\lam_ix_i\)\geqslant\sum_{i=1}^n\lam_i\psi(x_i).$$
In this article, we often take $p_1=\cdots=p_n=1/n$.
As special cases, if we consider the convex function $x^p$ with $p\in(-\8,0)\cup[1,+\8)$,
for $x_1,\cdots,x_n>0$,
$$\(\frac1n\sum_{i=1}^nx_i\)^p\leqslant\frac1n\sum_{i=1}^n x_i^p,\hs\mb{i.e.,}$$
\be\label{2.2A}
\(\sum_{i=1}^nx_i^p\)^{\frac1p}\geqslant n^{\frac1p-1}\sum_{i=1}^nx_i\mb{ for }p\geqslant1,\hs
\(\sum_{i=1}^nx_i^p\)^{\frac1p}\leqslant n^{\frac1p-1}\sum_{i=1}^nx_i\mb{ for }p<0;\ee
for the concave function $x^p$ with $p\in(0,1]$ and $x_1,\cdots,x_n>0$, we also have
\be\label{2.2B}\(\sum_{i=1}^nx_i^p\)^{\frac1p}
\leqslant n^{\frac1p-1}\sum_{i=1}^nx_i;\ee
for $\ln x$, which is a concave function, we have
\be\label{2.lala}\ln\sum_{i=1}^np_ix_i\geqslant\sum_{i=1}^np_i\ln x_i,\hs
\mb{i.e.,}\hs\sum_{i=1}^np_ix_i\geqslant\prod_{i=1}^nx_i^{p_i},\ee
where $x_1,\cdots,x_n$ are positive.
Actually \eqref{2.lala} can be regarded as a generalized version of mean value inequality.

We now recall the \textbf{Radon's inequality} in \cite{HL16,LZ11,L06,YKC85}
and their references, and it reads as follows:
if $a_i,b_i>0$, $i=1,\cdots,n$ and $m\in\R$, then
\be\sum_{i=1}^n\frac{a_i^{m+1}}{b_i^m}\geqslant
\frac{\disp\(\sum_{i=1}^na_i\)^{m+1}}{\disp\(\sum_{i=1}^nb_i\)^m},
\hs m\in(-\8,-1)\cup(0,+\8);\label{1.AA}\ee
\be\mb{and}\hs\sum_{i=1}^n\frac{a_i^{m+1}}{b_i^m}\leqslant
\frac{\disp\(\sum_{i=1}^na_i\)^{m+1}}{\disp\(\sum_{i=1}^nb_i\)^m},
\hs m\in(-1,0),\label{1.AB}\ee
where the equality ``$=$" only holds
when $\disp\frac{x_i}{y_i}=\cdots=\frac{x_n}{y_n}$.
Radon's inequality has been applied widely in high school education of mathematics
and International Mathematical Olympiads (IMO, see \cite{C13,HL16}).
Later, Radon's inequality was extended to the generalized form as follows.
\bt\label{GWSIP}
Let $a_i,b_i>0$, $i=1,\cdots,n$ and $p,q\in\R$.
If $q\in(-\8,-1)\cup[0,+\8)$, $p\geqslant q+1$ and $p(q+1)>0$, then
\be\sum_{i=1}^n\frac{a_i^p}{b_i^q}\geqslant n^{q+1-p}\cdot
\frac{\disp\(\sum_{i=1}^na_i\)^p}{\disp\(\sum_{i=1}^nb_i\)^q};\label{1.BA}\ee
if $q\in(-1,0]$ and $p\in(0,q+1]$, then
\be\sum_{i=1}^n\frac{a_i^p}{b_i^q}\leqslant n^{q+1-p}\cdot
\frac{\disp\(\sum_{i=1}^na_i\)^{p}}{\disp\(\sum_{i=1}^nb_i\)^q},
\label{1.BB}\ee
where the equality ``=" holds only when
\be\label{1.AC}\frac{x_i}{y_i}=\cdots=\frac{x_n}{y_n}.\ee\et
\bo For the reader's convenience, we provide a brief proof here.
We first consider the case when $q\in(-\8,-1)\cup[0,+\8)$, $p\geqslant q+1$ and $p(q+1)>0$.
We take $\tilde{a}_i=a_i^{\frac{p}{q+1}}$ and then by \eqref{1.AA},
$$\sum_{i=1}^n\frac{\tilde{a}_i^{q+1}}{b_i^q}\geqslant
\frac{\disp\left[\(\sum_{i=1}^na_i^{\frac{p}{q+1}}\)^{\frac{q+1}{p}}\right]^p}
{\disp\(\sum_{i=1}^nb_i\)^q}\geqslant
n^{q-p+1}\frac{\disp\(\sum_{i=1}^na_i\)^p}{\disp\(\sum_{i=1}^nb_i\)^q},$$
where the second ``$\geqslant$" follows by \eqref{2.2A} and \eqref{2.2B} from
$$\(\sum_{i=1}^na_i^{\frac{p}{q+1}}\)^{\frac{q+1}{p}}\geqslant
n^{\frac{q-p+1}{p}}\sum_{i=1}^na_i,\hs\mb{when }q\geqslant0,$$
$$\(\sum_{i=1}^na_i^{\frac{p}{q+1}}\)^{\frac{q+1}{p}}\leqslant
n^{\frac{q-p+1}{p}}\sum_{i=1}^na_i,\hs\mb{when }q<-1.$$
For the case when $q\in(-1,0]$ and $p\in(0,q+1]$, we similarly have
$$\sum_{i=1}^n\frac{\tilde{a}_i^{q+1}}{b_i^q}\leqslant
\frac{\disp\left[\(\sum_{i=1}^na_i^{\frac{p}{q+1}}\)^{\frac{q+1}{p}}\right]^p}
{\disp\(\sum_{i=1}^nb_i\)^q}\leqslant
n^{q-p+1}\frac{\disp\(\sum_{i=1}^na_i\)^p}{\disp\(\sum_{i=1}^nb_i\)^q},$$
where the second ``$\leqslant$" is obtained by \eqref{2.2B} and
$$\(\sum_{i=1}^na_i^{\frac{p}{q+1}}\)^{\frac{q+1}{p}}\leqslant
n^{\frac{q-p+1}{p}}\sum_{i=1}^na_i.$$
Thus the inequality \eqref{1.BB} is similarly obtained.
The proof is finished.
\eo

The following \textbf{Semiconcave-semiconvex Theorem} can be found
in \cite[Theorem 7.4]{H11}.
And this theorem is very effective to deduce more general inequalities.
\bt\label{th2.2}
Let $a<b$ and $x_1$, $\cdots$, $x_n\in[a,b]$ such that
\benu\item[(1)] $x_1\leqslant\cdots\leqslant x_n$;
\item[(2)] $x_1+\cdots+x_n=C$, where $C$ is a constant.\eenu
Let $f:[a,b]\ra\8$ be a function with $c\in(a,b)$ such that $f$ is concave (resp. convex)
on $[a,c]$ and convex (resp. concave) on $[c,b]$, and
$$F(x_1,\cdots,x_n)=f(x_1)+\cdots+f(x_n).$$
Then if $F$ achieves its minimum (resp. maximum) at some point $x=(x_1,\cdots,x_n)$,
then $x$ satisfies
$x_1=\cdots=x_{k-1}=a$, $x_{k+1}=\cdots=x_n$, $k=1,\cdots,n$;
if $F$ achieves its maximum (resp. minimum) at some point $x=(x_1,\cdots,x_n)$,
then $x$ satisfies
$x_1=\cdots=x_{k-1}$, $x_{k+1}=\cdots=x_n=b$, $k=1,\cdots,n$.
\et

In the sequel, when it comes to the derivatives of a function $f(x)$
on the bottom $a$ of an interval, we still use $f'(a)$ to
denote the unilateral derivatives for convenience if defined.
\section{Main inequalities}\label{s3}
\subsection{Main theorems}
In this subsection we are to present the main inequalities in various cases and prove them.
\bt\label{th3.1} Let $a_1,a_2,\cdots,a_n>0$, $m,p,t,r,\beta\in\R$ with $t\geqslant0$,
$p\ne0$ and
$$s=\sum_{i=1}^na_i^p\hs\mb{ with }\hs ts>ra_i^p\hs\mb{for each }i=1,\cdots,n.$$
Let
$$T=\left\{\ba{ll}s,&\mb{if }r\leqslant t,\\ ts/r,&\mb{if }r>t\ea\right.$$
and $f:\R\ra\R$ be a parabolic function such that
$$f(x)=x^2+\frac{2m}{(\beta+1)rp}x+\frac{m(m-p)}{\beta(\beta+1)r^2p^2},$$
with $\beta$ and $r$ chosen appropriately.
When there exist two different real solutions to the equation $g(x)=0$,
we set $X_1$ and $X_2$ to be the two solutions with $X_1<X_2$, i.e.,
\be\label{3.A}X_1=-\frac{m}{(\beta+1)rp}-\sqrt{\frac{m(\beta p+p-m)}{\beta(\beta+1)^2r^2p^2}}
\hs\mb{and}\hs
X_2=-\frac{m}{(\beta+1)rp}+\sqrt{\frac{m(\beta p+p-m)}{\beta(\beta+1)^2r^2p^2}}.\ee

Then we have the following conclusions.

\noindent(1) Suppose that $\beta$, $r$, $m$, $p$ and $t\in\R$ satisfy
each one of the following four cases:
\benu\item[(i)] $\beta r=0$ and $m(m-p)\geqslant0$;
\item[(ii)] $\beta=-1$, $r\ne0$ and either
\benu\item[(ii.1)] $m(m-p)t\geqslant0$ and $m(m+p)r\leqslant0$, or
\item[(ii.2)] $m(m-p)t>0$, $m(m+p)r>0$ and $\disp\frac{(m-p)ts}{(m+p)r}\geqslant T$;\eenu
\item[(iii)] $\beta\in(-\8,-1)\cup(0,+\8)$, $r\ne0$ and one of the following cases holds,
\benu\item[(iii.1)] $\beta pm>0$ and $p[m-(\beta+1)p]\geqslant0$,
\item[(iii.2)] $\beta m[m-(\beta+1)p]\geqslant0$,
\item[(iii.3)] $mrp\beta\geqslant0$ and $m(m-p)\geqslant 0$,
\item[(iii.4)] $t>r$, $mrp\beta<0$, $m(m-p)>0$, $\beta m[m-(\beta+1)p]<0$
and $(t-r)X_1\geqslant1$;
\eenu
\item[(iv)] $\beta\in(-1,0)$, $t>r\ne0$, $m(m-p)\geqslant0$, $m[m-(\beta+1)p]>0$
and $(t-r)X_2\geqslant1$.\eenu
Then
\be\sum_{i=1}^n\frac{a_i^m}{(ts-ra_i^p)^\beta}\geqslant
\frac{n^{\beta+1-\frac mp}}{(nt-r)^{\beta}}\(\sum_{i=1}^na_i^p\)^{\frac mp-\beta}.
\label{3.zz1}\ee

\noindent(2) Suppose that $\beta$, $r$, $m$, $p$ and
$t\in\R$ satisfy one of the following cases:
\benu\item[(v)] $\beta r=0$ and $m(m-p)\leqslant0$;
\item[(vi)] $\beta=-1$, $r\ne0$ and either
\benu\item[(vi.1)] $m(m-p)t\leqslant0$ and $m(m+p)r\geqslant0$, or
\item[(vi.2)] $m(m-p)t<0$, $m(m+p)r<0$ and $\disp\frac{(m-p)ts}{(m+p)r}\geqslant T$;\eenu
\item[(vii)] $\beta\in(-1,0)$ and one of the following cases holds,
\benu\item[(vii.1)] $0\leqslant pm\leqslant(\beta+1)p^2$,
\item[(vii.2)] $mrp\geqslant0$ and $m(m-p)\geqslant 0$,
\item[(vii.3)] $r<0$, $(\beta+1)p^2<pm\leqslant p^2$ and $(t-r)X_1\geqslant1$;
\eenu
\item[(viii)] $\beta\in(-\8,-1)\cup(0,+\8)$, $t>r\ne0$, $m(m-p)\leqslant0$,
$\beta m[m-(\beta+1)p]<0$ and $(t-r)X_2\geqslant1$.
\eenu
Then
\be\sum_{i=1}^n\frac{a_i^m}{(ts-ra_i^p)^\beta}\leqslant
\frac{n^{\beta+1-\frac mp}}{(nt-r)^{\beta}}
\(\sum_{i=1}^na_i^p\)^{\frac mp-\beta}.\label{3.zz2}\ee
\et
\bo Noting that
$$ts>r\max_{1\leqslant i\leqslant n}\{a_i^p\},$$
we know that
$$nt\geqslant\frac{ts}{\disp\max_{1\leqslant i\leqslant n}\{a_i^p\}}>r,$$
which implies that the right sides of \eqref{3.zz1} and \eqref{3.zz2} make sense.

We consider the function $g:(0,T)\ra\R$ with  such that
\be\label{3.yy1}g(x)=\frac{x^{\frac mp}}{(ts-rx)^\beta}.\ee
Then we know that for each $x\in(0,T)$,
\be\label{3.yy2}g''(x)=\frac{x^{\frac mp-2}}{(ts-rx)^{\beta}}\left[\frac mp\(\frac mp-1\)
+\frac{2m\beta r}p\frac{x}{ts-rx}+\beta(\beta+1)r^2\frac{x^2}{(ts-rx)^2}\right].\ee
If $\beta=-1$,
\be\label{3.yy3}g''(x)=\frac{m}{p^2}x^{\frac mp-2}\left[(m-p)ts
-(m+p)rx\right];\ee
if $\beta r=0$,
\be\label{3.3}g''(x)=\frac mp\(\frac mp-1\)\frac{x^{\frac mp-2}}{(ts-rx)^{\beta}};\ee
if $\beta(\beta+1)r\ne0$,
\begin{align}\notag g''(x)=&(\beta+1)\beta r^2\frac{x^{\frac mp-2}}{(ts-rx)^{\beta}}
f\(\frac{x}{ts-rx}\)\\
=&(\beta+1)\beta r^2\frac{x^{\frac mp-2}}{(ts-rx)^{\beta}}
\left[\(\frac{x}{ts-rx}+\frac{m}{(\beta+1)rp}\)^2
+\frac{m(m-\beta p-p)}{\beta(\beta+1)^2r^2p^2}\right].\label{3.4}
\end{align}

In the following, we divide it into three parts to show the conclusions.

\noindent\textbf{Part 1.} We first show the conclusions for the cases (i), (ii),
(iii.2), (iii.3), (v), (vi), (vii.1) and (vii.2).

According to \eqref{3.yy3}, \eqref{3.3} and \eqref{3.4},
we know that when each one of the cases (i), (ii)
and the case when $\beta(\beta+1)>0$, $r\ne0$,
\be\label{3.5}\mb{and either}\hs\frac{m[m-(\beta+1)p]}{\beta}\geqslant0,\ee
\be\label{3.6}\hs\mb{or}\hs\frac{m}{(\beta+1)rp}\geqslant0\hs\mb{and}\hs
\frac{m(m-p)}{\beta(\beta+1)r^2p^2}\geqslant0\ee
hold, $g(x)$ is a convex function on $(0,T)$;
and when each one of the cases (v), (vi) and the case when $\beta\in(-1,0)$, $r\ne0$ and
\eqref{3.5} hold, $g(x)$ is a concave function on $(0,T)$.
Hence, when (i), (ii), \eqref{3.5} or \eqref{3.6} holds,
we employ the Jensen's inequality and obtain
$$\sum_{i=1}^n\frac{a_i^m}{(ts-ra_i^p)^\beta}\geqslant
n\frac{(s/n)^{\frac mp}}{(ts-rs/n)^\beta}
=\frac{n^{\beta+1-\frac mp}}{(nt-r)^{\beta}}\(\sum_{i=1}^na_i^p\)^{\frac mp-\beta},$$
which is exactly \eqref{3.zz1} and when (v), (vi) or \eqref{3.5} holds,
we similarly have \eqref{3.zz2}.
Noting that in case when $\beta(\beta+1)>0$,
\eqref{3.5} is equivalent to (iii.2),
\eqref{3.6} is equivalent to (iii.3)
and in case when $\beta\in(-1,0)$,
\eqref{3.5} is equivalent to (vii.1),
\eqref{3.6} is equivalent to (vii.2),
we can see that the conclusions for the cases (iii.2),
(iii.3), (vii.1) and (vii.2) with $r\ne0$ have been proved.
\vs

\noindent\textbf{Part 2.} Next, we show the conclusions for (iii.4), (iv) and (viii).
Set $t>r$ and $\beta(\beta+1)r\ne0$.
In consideration of the parabolic function $f(x)$,
there are obviously some other cases such that $g''(x)\geqslant0$ on $(0,T)$ by adjusting
the axis of symmetry for $f(x)$, $y$-intercept of $f(x)$ and the solutions $X_1$, $X_2$:
\be\label{3.10}\beta(\beta+1)>0,\hs\frac{m}{(\beta+1)rp}<0,\hs
\frac{m(m-p)}{\beta(\beta+1)r^2p^2}>0\hs\mb{and}\hs X_1\geqslant\frac1{t-r};\ee
\be\label{3.11}\beta\in(-1,0),\hs\frac{m(m-p)}{\beta(\beta+1)r^2p^2}\leqslant0\hs\mb{and}
\hs X_2\geqslant\frac1{t-r}.\ee
For the existence of $X_1$ and $X_2$, it is also required that
\be\label{3.12}\beta m[m-(\beta+1)p]<0.\ee
Noting that \eqref{3.10} and \eqref{3.12} $\Lra$ (iii.4) and
\eqref{3.11} and \eqref{3.12} $\Lra$ (iv), we can similarly obtain
\eqref{3.zz1} for the cases (iii.4) and (iv).
The situation for the cases (vii.3) and (viii) can be similarly guaranteed.
\vs

\noindent\textbf{Part 3.} At last, it remains to prove the conclusion for (iii.1).
Indeed, in this case $\disp\frac mp\geqslant\beta+1$ and $(\beta+1)pm>0$.
Then by generalized Radon's inequality (Theorem \ref{GWSIP}), one sees
$$\sum_{i=1}^n\frac{a_i^m}{(ts-ra_i^p)^\beta}=\sum_{i=1}^n\frac{(a_i^p)^{\frac mp}}
{(ts-ra_i^p)^\beta}\geqslant\frac{\disp n^{\beta-\frac mp+1}\(\sum_{i=1}^na_i^p\)^{\frac mp}}
{\disp\(\sum_{i=1}^n(ts-ra_i^p)\)^{\beta}}=\frac{n^{\beta+1-\frac mp}}
{(nt-r)^\beta}\(\sum_{i=1}^na_i^p\)^{\frac mp-\beta}.$$
The proof is hence accomplished now.
\eo

\br In Theorem \ref{th3.1}, some different cases have non-empty intersections,
but for writing brevity, we do not classify them explicitly.
\er
Next under the conditions of Theorem \ref{th3.1}, we compare
$$\sum_{i=1}^n\frac{a_i^m}{(ts-ra_i^p)^\beta}\hs\mb{and}\hs
\frac1{(nt-r)^\beta}\sum_{i=1}^na_i^{m-\beta p}.$$
First, we observe that when $\beta=0$, $p=0$ or $t=0$,
it always holds that
$$\sum_{i=1}^n\frac{a_i^m}{(ts-ra_i^p)^\beta}=
\frac1{(nt-r)^\beta}\sum_{i=1}^na_i^{m-\beta p}.$$
Hence we only consider the cases when $\beta pt\ne0$ in the following.

\bt\label{th3.2}
Under the conditions of Theorem \ref{th3.1} with $\beta pt\ne0$,
we have the inequality
\be\sum_{i=1}^n\frac{a_i^m}{(ts-ra_i^p)^\beta}\geqslant
\frac1{(nt-r)^\beta}\sum_{i=1}^na_i^{m-\beta p}\label{3.1}\ee
in the following cases:
\benu\item[(i)] $\beta>0$ and one of the following cases holds,
\benu\item[(i.1)] $\beta\geqslant 1$, $r\geqslant0$ and $\beta p^2\leqslant pm$
or $m=(\beta+1)p$,
\item[(i.2)] $\beta\geqslant 1$, $r<0$ and $(\beta+1)p^2\leqslant pm$,
\item[(i.3)] $\beta\in(0,1)$, $r\geqslant0$ and $p^2\leqslant pm\leqslant(\beta+1)p^2$,
\item[(i.4)] $r<0$, $(\beta\vee1)p^2\leqslant pm<(\beta+1)p^2$ and $(t-r)X_1\geqslant1$;
\eenu
\item[(ii)] $\beta\in(-1,0)$ and one of the following cases holds,
\benu\item[(ii.1)] $r=0$ and $\beta p^2<pm\leqslant0$,
\item[(ii.2)] $t\geqslant r$ and $pm\leqslant\beta p^2$,
\item[(ii.3)] $t>r\ne0$, $\beta p^2<pm\leqslant0$ and $(t-r)X_2\geqslant 1$;
\eenu
\item[(iii)] $\beta=-1$ and one of the following cases holds,
\benu\item[(iii.1)] $pm\leqslant-p^2$, or $m=0$,
\item[(iii.2)] $r\geqslant0$ and $-p^2<pm<0$,
\item[(iii.3)] $r<0$, $-p^2<pm<0$ and $\disp\frac{(m-p)t}{(m+p)r}\geqslant1$;
\eenu
\item[(iv)] $\beta<-1$ and one of the following cases holds,
\benu\item[(iv.1)] $pm\leqslant\beta p^2$, or $m=(\beta+1)p$,
\item[(iv.2)] $r\geqslant0$ and $\beta p^2<pm<(\beta+1)p^2$,
\item[(iv.3)] $r<0$, $\beta p^2<pm<(\beta+1) p^2$ and $(t-r)X_1\geqslant1$,\eenu\eenu
where $a\vee b$ means the bigger one of $a,b\in\R$ and $X_1$ and $X_2$ are given in \eqref{3.A}.
\et
\bo We first show \eqref{3.1} in the cases (i.1) (with $\beta p^2\leqslant pm$), (i.2),
(ii.2), (iii.1) and (iv.1) with
$pm\leqslant\beta p^2$ in the first three parts and for other cases in the fourth part.

\noindent\textbf{Part 1.} (1) We first consider (i.1) (with $\beta p^2\leqslant pm$)
and (i.2) in this part and prove \eqref{3.1} in the following cases in advance:
\benu\item[(i.1a)] $\beta\geqslant1$, $r\geqslant0$, $p>0$ and $\beta p\leqslant m$;
\item[(i.2a)] $\beta\geqslant1$, $r<0$, $p>0$ and $(\beta+1)p\leqslant m$.
\eenu

We can assume that $a_1\leqslant a_2\leqslant\cdots\leqslant a_n$
for writing convenience.
Set
\be\label{2.llm}A_i=\frac{a_i^{m-p}}{(ts-ra_i^p)^{\beta}}.\ee
Then we know that $a_1^p\leqslant a_2^p\leqslant\cdots\leqslant a_n^p$ and
$A_1\leqslant A_2\leqslant\cdots\leqslant A_n$,
since the function $x^{m-p}/(ts-rx^p)^\beta$ is non-decreasing in $x>0$ in each case.
By the rearrangement inequality, we have
\be t\sum_{i=1}^na_i^pA_i\geqslant t\sum_{i=1}^na_{i+k}^pA_i,\hs\mb{for all }k=1,\cdots,n-1,
\label{2.lll}\ee
\be\mb{and}\hs (t-r)\sum_{i=1}^na_i^pA_i=(t-r)\sum_{i=1}^na_{i}^pA_i,\label{2.llk}\ee
where when $i+k>n$, $a_{i+k}$ is taken to be $a_{i+k-n}$.
Adding all inequalities in \eqref{2.lll} for each $k=1,\cdots,n-1$ and \eqref{2.llk} up,
we obtain
$$(nt-r)\sum_{i=1}^na_i^pA_i\geqslant\sum_{i=1}^n(ts-ra_i^p)A_i$$
\be\label{2.llj}\mb{and}\hs\sum_{i=1}^n\frac{a_i^m}{(ts-ra_i^p)^\beta}\geqslant
\frac1{nt-r}\sum_{i=1}^n\frac{a_i^{m-p}}{(ts-ra_i^p)^{\beta-1}}.
\ee

Now when $\beta\geqslant1$, noticing that for each $r\in\R$,
\be\label{3.kkc}a_1^{m-\beta p}\leqslant a_2^{m-\beta p}\leqslant\cdots
\leqslant a_n^{m-\beta p}\hs\mb{and}\ee
\be\label{3.kkd}\(\frac{a_1^p}{ts-ra_1^p}\)^\beta
\leqslant\(\frac{a_2^p}{ts-ra_2^p}\)^\beta\leqslant\cdots
\leqslant\(\frac{a_n^p}{ts-ra_n^p}\)^\beta,\ee
we can adopt the Chepyshev's inequality and have
\be\label{2.llh}\sum_{i=1}^n\frac{a_i^m}{(ts-ra_i^p)^\beta}
=\sum_{i=1}^na_i^{m-\beta p}\(\frac{a_i^p}{ts-ra_i^p}\)^\beta
\geqslant\sum_{i=1}^na_i^{m-\beta p}\cdot
\frac1n\sum_{i=1}^n\(\frac{a_i^p}{ts-ra_i^p}\)^\beta.
\ee
Since $x^\beta$ is a convex increasing function,
we can use the Jensen's inequality and obtain that
\be\label{2.llg}\frac1n\sum_{i=1}^n\(\frac{a_i^p}{ts-ra_i^p}\)^\beta
\geqslant\(\frac1n\sum_{i=1}^n\frac{a_i^p}{ts-ra_i^p}\)^\beta
\geqslant\frac{1}{(nt-r)^\beta},\ee
where we have used \eqref{2.llj} by setting $\beta=1$ and $m=p$.
We have actually obtained \eqref{3.1} for these two cases
by combining \eqref{2.llh} and \eqref{2.llg}.
\vs

\noindent\textbf{(2)} We then prove \eqref{3.1} for the cases:
\benu\item[(i.2b)] $\beta\geqslant1$, $r\geqslant0$, $p<0$ and $m\leqslant\beta p$;
\item[(ii.2b)] $\beta\geqslant1$, $r<0$, $p<0$ and $m\leqslant(\beta+1)p$.
\eenu

We also assume that $a_1\leqslant a_2\leqslant\cdots\leqslant a_n$ and
set $A_i$ as \eqref{2.llm}.
Then we know
\be\label{2.llak}a_1^p\geqslant a_2^p\geqslant\cdots\geqslant a_n^p\hs\mb{and}
\hs A_1\geqslant A_2\geqslant\cdots\geqslant A_n,\ee
since the functions $x^p$ and $x^{m-p}/(ts-rx^p)^\beta$ are both non-increasing
in $x>0$ in each case.
By the rearrangement inequality, we can similarly obtain \eqref{2.lll},
\eqref{2.llk} and \eqref{2.llj}.

Then if $\beta\in\N^+$, we can similarly deduce \eqref{2.llj}.
If $\beta\geqslant1$, noticing also that
\be\label{2.llc}a_1^{m-\beta p}\geqslant a_2^{m-\beta p}\geqslant\cdots\geqslant
a_n^{m-\beta p}\hs\mb{and}\ee
\be\label{2.llb}\(\frac{a_1^p}{ts-ra_1^p}\)^\beta\geqslant
\(\frac{a_2^p}{ts-ra_2^p}\)^\beta\geqslant\cdots
\geqslant\(\frac{a_n^p}{ts-ra_n^p}\)^\beta,\ee
we can also adopt the Chepyshev's inequality and have \eqref{2.llh}, \eqref{2.llg}
and then \eqref{3.1}, finally.
\vs

\noindent\textbf{Part 2.} Now we consider the case (ii.2).
We first prove \eqref{3.1} for the case when $\beta\in(-1,0)$, $t\geqslant r$,
$p>0$ and $m\leqslant\beta p$.
Similarly, we assume that $a_1\leqslant a_2\leqslant\cdots\leqslant a_n$.
Noticing that \eqref{2.llc} and \eqref{2.llb} still hold for this case,
and then we again obtain \eqref{2.llh}.
By the mean value inequality, we can see that
\be\label{2.kke}\frac1n\sum_{i=1}^n\(\frac{a_i^p}{ts-ra_i^p}\)^\beta
\geqslant\(\prod_{i=1}^n\frac{a_i^p}{ts-ra_i^p}\)^{\frac{\beta}n}
=\(\frac{\disp\prod_{i=1}^n a_i^p}
{\disp\prod_{i=1}^n(ts-ra_i^p)}\)^{\frac{\beta}n}\ee
and by \eqref{2.lala},
\begin{align*}ts-ra_i^p=&ta_1^p+\cdots+ta_{i-1}^p+(t-r)a_i^p+ta_{i+1}^p+\cdots+ta_n^p\\
\geqslant&(nt-r)a_i^{\frac{p(t-r)}{nt-r}}\prod_{k=1,k\ne i}^na_k^{\frac{pt}{nt-r}}.
\end{align*}
Hence by\eqref{2.kke}
\be\label{3.kkb}\frac1n\sum_{i=1}^n\(\frac{a_i^p}{ts-ra_i^p}\)^\beta
\geqslant\(\frac{\disp\prod_{i=1}^n a_i^p}
{\disp(nt-r)^n\prod_{k=1}^na_k^{\frac{p(nt-r)}{nt-r}}}\)^{\frac{\beta}n}=(nt-r)^{-\beta}.\ee
and then \eqref{3.1} follows from \eqref{3.kkb} and \eqref{2.llh}.

For the case when $\beta\in(-1,0)$, $t\geqslant r$, $p<0$ and $m\geqslant\beta p$,
we can similarly deduce \eqref{3.kkc}, \eqref{3.kkd} and \eqref{2.llh}.
Then since $\beta<0$, \eqref{2.kke} and \eqref{3.kkb} are also valid
and \eqref{3.1} holds true for this case.

\noindent\textbf{Part 3.} We then consider the cases (iii.1) and (iv.1) with $pm\leqslant\beta p^2$.
We first consider the case when $\beta\leqslant-1$, $r\in\R$, $p>0$ and $m\leqslant\beta p$.
We similarly assume that $a_1\leqslant a_2\leqslant\cdots\leqslant a_n$
and then for all $r\in\R$, \eqref{2.llc}, \eqref{2.llh} and \eqref{2.llb} are valid.
Since $x^{-\beta}$ is a convex increasing function,
we can use the Jensen's inequality and obtain that
\begin{align}\notag\frac1n\sum_{i=1}^n\(\frac{a_i^p}{ts-ra_i^p}\)^\beta=&
\frac1n\sum_{i=1}^n\(t\sum_{k=1}^n\frac{a_k^p}{a_i^p}-r\)^{-\beta}
\geqslant\(\frac tn\sum_{k=1}^n\sum_{i=1}^n\frac{a_k^p}{a_i^p}-r\)^{-\beta}\\
\geqslant&\(\frac{t}nn^2-r\)^{-\beta}=\frac1{(nt-r)^{\beta}}.\label{2.ller}\end{align}
Then \eqref{3.1} is proved for all $\beta\leqslant-1$
by combining \eqref{2.llh} and \eqref{2.ller}.

For the case when $\beta\leqslant-1$, $r\in\R$, $p<0$ and $m\geqslant\beta p$,
we can similarly obtain \eqref{3.kkc} and \eqref{3.kkd}.
Then with the same argument as above, we can show \eqref{3.1} in this case.
\vs

\noindent\textbf{Part 4.} We now consider other cases,
in which cases we adopt Theorem \ref{th3.1} to show \eqref{3.1}.
Actually, in other cases, it holds that $(m-\beta p)/p\in[0,1]$.
This implies that the function $x^{\frac{m}{p}-\beta}$ is concave.
Then by Jensen's inequality, we deduce that
\be\label{3.RRR}\(\sum_{i=1}^na_i^p\)^{\frac mp-\beta}
\geqslant n^{\frac{m}{p}-1-\beta}\sum_{i=1}^na_i^{m-\beta p}.\ee
Next we observe that (i.1) with $m=(\beta+1)p$ satisfies (iii.1) of Theorem \ref{th3.1},
(i.3) and (iv.2) satisfies (i) or (iii.3) of Theorem \ref{th3.1},
(i.4) and (iv.3) satisfy the case (iii.4) of Theorem \ref{th3.1},
(ii.1) satisfies (i) of Theorem \ref{th3.1},
(ii.3) satisfies (iv) of Theorem \ref{th3.1},
(iii.1) (with $m=0$) and (iii.2) satisfy (ii.1) of Theorem \ref{th3.1},
(iii.3) satisfies (ii.2) of Theorem \ref{th3.1},
(iv.1) (with $m=(\beta+1)p$) satisfies (iii.2) of Theorem \ref{th3.1},
and (iv.3) satisfies (iii.4) of Theorem \ref{th3.1}.
Then \eqref{3.1} follows from \eqref{3.RRR} and the result \eqref{3.zz1} in these cases.
The proof is complete.
\eo

\bt\label{th3.3} Under the conditions of Theorem \ref{th3.1} with $\beta pt\ne0$,
we have the inequality
\be\sum_{i=1}^n\frac{a_i^m}{(ts-ra_i^p)^\beta}\leqslant
\frac1{(nt-r)^\beta}\sum_{i=1}^na_i^{m-\beta p},\label{3.2}\ee
in the following cases:
\benu\item[(i)] $\beta>0$ and either
\benu\item[(i.1)] $r\leqslant0$ and $pm\leqslant(\beta\wedge1)p^2$, or
\item[(i.2)] $t>r>0$, $0\leqslant pm\leqslant(\beta\wedge 1)p^2$ and $(t-r)X_2\geqslant1$;
\eenu
\item[(ii)] $\beta\in(-1,0)$ and one of the following cases holds,
\benu\item[(ii.1)] $m=(\beta+1)p$,
\item[(ii.2)] $r>0$ and $pm\geqslant p^2$,
\item[(ii.3)] $r<0$ and $pm\leqslant\beta p^2$,
\item[(ii.4)] $r<0$, $(\beta+1)p^2<pm\leqslant p^2$ and $(t-r)X_1\geqslant1$;\eenu
\item[(iii)] $\beta\leqslant-1$ and one of the following cases holds,
\benu\item[(iii.1)] $\beta=-1$, $r>0$ and $0\leqslant pm\leqslant p^2$,
\item[(iii.2)] $r=0$ and $0\leqslant pm\leqslant p^2$,
\item[(iii.3)] $\beta\in\Z\sm\N$, $r<0$ and $pm\geqslant0$,
\item[(iii.4)] $\beta<-1$, $t>r>0$, $0\leqslant pm\leqslant p^2$ and $(t-r)X_2\geqslant1$,
\eenu\eenu
where $a\wedge b$ means the smaller one of $a,b\in\R$
and $X_2$ is given in \eqref{3.A}.
\et
\bo We split the proof into three parts.

\noindent\textbf{Part 1.} We first show \eqref{3.2} for (i.1) and the case
when $\beta\geqslant1$, $r\leqslant 0$, $p>0$ and $p\geqslant m$.
In this case, the function $x^p$ is non-decreasing and $x^{m-p}/(ts-rx^p)^\beta$
is non-increasing.
Hence by setting $a_1\leqslant a_2\leqslant\cdots\leqslant a_n$ and \eqref{2.llm},
we have
\be\label{2.llai}a_1^p\leqslant a_2^p\leqslant\cdots\leqslant a_n^p\hs\mb{and}
\hs A_1\geqslant A_2\geqslant\cdots\geqslant A_n.\ee
By the rearrangement inequality, we obtain
\be t\sum_{i=1}^na_i^pA_i\leqslant t\sum_{i=1}^na_{i+k}^pA_i,\hs\mb{for all }k=1,\cdots,n-1,
\label{2.llah}\ee
\be\mb{and}\hs (t-r)\sum_{i=1}^na_i^pA_i=(t-r)\sum_{i=1}^na_{i}^pA_i,\label{2.llag}\ee
where it is also taken that $a_{i+k}=a_{i+k-n}$, when $i+k>n$.
Adding all inequalities in \eqref{2.llah} for each $k=1,\cdots,n-1$ and \eqref{2.llag} up,
we can further obtain
\be\label{2.llaf}\sum_{i=1}^n\frac{a_i^m}{(ts-ra_i^p)^\beta}\leqslant
\frac1{nt-r}\sum_{i=1}^n\frac{a_i^{m-p}}{(ts-ra_i^p)^{\beta-1}}.
\ee

Next we consider $\beta\in[0,1)$, $r\leqslant0$, $p>0$ and $\beta p\geqslant m$.
Since $x^{m-\beta p}$ is non-increasing and $x/(ts-rx)$ is non-decreasing,
we have
\be\label{3.mma}a_1^{m-\beta p}\geqslant a_2^{m-\beta p}\geqslant\cdots
\geqslant a_n^{m-\beta p}\hs\mb{and}\ee
\be\label{3.mmb}\(\frac{a_1^p}{ts-ra_1^p}\)^\beta
\leqslant\(\frac{a_2^p}{ts-ra_2^p}\)^\beta\leqslant\cdots
\leqslant\(\frac{a_n^p}{ts-ra_n^p}\)^\beta.\ee
By the Chepyshev's inequality, we obtain
\be\label{3.mmc}\sum_{i=1}^n\frac{a_i^m}{(ts-ra_i^p)^\beta}
\leqslant\sum_{i=1}^na_i^{m-\beta p}\cdot
\frac1n\sum_{i=1}^n\(\frac{a_i^p}{ts-ra_i^p}\)^\beta.
\ee
By Jensen's inequality, we have
\be\label{3.mmd}\frac1n\sum_{i=1}^n\(\frac{a_i^p}{ts-ra_i^p}\)^\beta\leqslant
\(\frac1n\sum_{i=1}^n\frac{a_i^p}{ts-ra_i^p}\)^{\beta}
\leqslant\frac{1}{(nt-r)^{\beta}},\ee
where we used \eqref{2.llaf} with $\beta=1$ and $m=p$.
Then \eqref{3.2} follows from \eqref{3.mmc} and \eqref{3.mmd}.

Now we consider $\beta\geqslant1$, $r\leqslant0$, $p>0$ and $p\geqslant m$.
One can see by \eqref{2.llaf} that
\be\label{3.mme}\sum_{i=1}^n\frac{a_i^m}{(ts-ra_i^p)^\beta}\leqslant
\frac1{nt-r}\sum_{i=1}^n\frac{a_i^{m-p}}{(ts-ra_i^p)^{\beta-1}}\leqslant\cdots
\leqslant\frac1{(nt-r)^{\lfloor\beta\rfloor}}
\sum_{i=1}^n\frac{a_i^{m-\lfloor\beta\rfloor p}}{(ts-ra_i^p)^{\beta-\lfloor\beta\rfloor}},\ee
where $\lfloor\beta\rfloor$ is the largest integer no more than $\beta$.
Noting that $0\leqslant\beta-\lfloor\beta\rfloor<1$,
we infer from \eqref{3.2} for the case (i.1) with $\beta\in[0,1)$ that
\be\label{3.mmf}\frac1{(nt-r)^{\lfloor\beta\rfloor}}
\sum_{i=1}^n\frac{a_i^{m-\lfloor\beta\rfloor p}}{(ts-ra_i^p)^{\beta-\lfloor\beta\rfloor}}
\leqslant\frac1{(nt-r)^{\beta}}
\sum_{i=1}^n a_i^{m-\beta p},\ee
which is exactly \eqref{3.2} in this case.
\vs

Hereafter we consider the case when $\beta>0$, $r\leqslant 0$, $p<0$
and $\beta p\leqslant m$ and prove \eqref{3.2} for $\beta\geqslant1$ first.
In this case the function $x^p$ is non-increasing and $x^{m-p}/(ts-rx^p)^\beta$ is non-decreasing.
Similarly by setting $a_1\leqslant a_2\leqslant\cdots\leqslant a_n$ and \eqref{2.llm},
\be\label{3.42}a_1^p\geqslant a_2^p\geqslant\cdots\geqslant a_n^p\hs\mb{and}
\hs A_1\leqslant A_2\leqslant\cdots\leqslant A_n,\ee
which implies \eqref{2.llaf} in this case.
Then for $\beta\in[0,1)$, $r\leqslant0$, $p<0$ and $\beta p\leqslant m$,
we have
\be\label{3.43}a_1^{m-\beta p}\leqslant a_2^{m-\beta p}\leqslant\cdots
\leqslant a_n^{m-\beta p}\hs\mb{and}\ee
\be\label{3.44}\(\frac{a_1^p}{ts-ra_1^p}\)^\beta
\geqslant\(\frac{a_2^p}{ts-ra_2^p}\)^\beta\geqslant\cdots
\geqslant\(\frac{a_n^p}{ts-ra_n^p}\)^\beta,\ee
and \eqref{3.mmd} and \eqref{3.mme} are obtained.
Hence \eqref{3.2} is proved.
For $\beta\geqslant1$, $r\leqslant0$, $p<0$ and $p\leqslant m$, \eqref{3.2}
can be deduced by similar argument.
\vs

\noindent\textbf{Part 2.} Next we prove \eqref{3.2} for other cases except (iii.3),
each of which satisfies
\be\label{3.45}\frac{m-\beta p}{p}\in(-\8,0]\cup[1,+\8).\ee

For these cases, we need to employ the results from Theorem \ref{th3.1}.
Note that as long as these cases satisfy the cases in (2) of Theorem \ref{th3.1}
and \eqref{3.45}, which implies
\be\label{3.46}\(\sum_{i=1}^na_i^p\)^{\frac mp-\beta}\leqslant
n^{\frac{m}{p}-1-\beta}\sum_{i=1}^na_i^{m-\beta p},\ee
the inequality \eqref{3.2} follows directly from \eqref{3.46}.

Actually, it is not hard to check that (i.2) and (iii.4) satisfy (viii) of Theorem \ref{th3.1},
(ii.1) satisfies (vii.1) of Theorem \ref{th3.1},
(ii.2) and (ii.3) satisfy (vii.2) of Theorem \ref{th3.1},
(ii.4) satisfies (vii.3) of Theorem \ref{th3.1},
(iii.1) satisfies (vi.1) of Theorem \ref{th3.1}
and (iii.2) satisfies (v) of Theorem \ref{th3.1}.
The proof of Part 2 ends here thereby.
\vs

\noindent\textbf{Part 3.} At last we show \eqref{3.2} in the case (iii.3).
We consider the case when $\beta\leqslant-1$, $r<0$ and $m,p\geqslant0$ in advance.
Indeed, in this case, we can also see that the functions $x^p$ and $x^{m}/(ts-rx^p)^{\beta+1}$
are both non-decreasing in $x>0$.
Then similar to the argument from \eqref{2.llm} to \eqref{2.llj},
we have
\be\label{3.47}a_1^p\leqslant a_2^p\leqslant\cdots\leqslant a_n^p\hs\mb{and}\ee
\be\label{3.48}\frac{a_1^m}{(ts-ra_1^p)^{\beta+1}}\leqslant\frac{a_2^m}{(ts-ra_2^p)^{\beta+1}}
\leqslant\cdots\leqslant\frac{a_n^m}{(ts-ra_n^p)^{\beta+1}}.\ee
Analogously, we obtain
$$(nt-r)\sum_{i=1}^n\frac{a_i^{m+p}}{(ts-ra_i^p)^{\beta+1}}\geqslant
\sum_{i=1}^n\frac{a_i^m}{(ts-ra_i^p)^{\beta}},$$
\be\label{3.49}\mb{i.e.,}\hs\sum_{i=1}^n\frac{a_i^{m}}{(ts-ra_i^p)^{\beta}}\leqslant
\frac{1}{(nt-r)^{-1}}\sum_{i=1}^n\frac{a_i^{m+p}}{(ts-ra_i^p)^{\beta+1}}.\ee

Then if $\beta\leqslant-1$,
\be\label{3.50}\sum_{i=1}^n\frac{a_i^{m}}{(ts-ra_i^p)^{\beta}}\leqslant\cdots\leqslant
\frac{1}{(nt-r)^{\lceil\beta\rceil}}
\sum_{i=1}^n\frac{a_i^{m-\lceil\beta\rceil p}}{(ts-ra_i^p)^{\beta-\lceil\beta\rceil}},\ee
where $\lceil\beta\rceil$ means the smallest integer no less than $\beta$.
Therefore, when $\beta\in\Z\sm\N$, \eqref{3.50} is exactly \eqref{3.2}.

For the case when $\beta\in\Z\sm\N$, $r<0$ and $m,p\leqslant0$,
the function $x^p$ and $x^{m}/(ts-rx^p)^{\beta+1}$ are both non-increasing in $x>0$.
Hence \eqref{3.47} and \eqref{3.48} hold with all $\leqslant$'s replaced by $\geqslant$'s.
Then \eqref{3.49} and \eqref{3.50} are also valid.
Finally, we can similarly deduce \eqref{3.2} in this case.
The proof is complete now.\eo

\br In the proofs of Theorems \ref{th3.2} and \ref{th3.3},
we also obtain some interesting inequalities \eqref{2.llj} and \eqref{3.49}
in corresponding cases presented therein.
These inequalities can also be viewed as generalizations of those obtained in \cite{JG23}.
\er

\subsection{Other cases}

Up to now we have proved the main theorems, but there are still other cases
which can guarantee the inequalities \eqref{3.zz1}, \eqref{3.zz2}, \eqref{3.1} and \eqref{3.2}.
For example, when $t>r\ne0$ and one of the followings holds:
\benu\item[(A)] $\beta\in(-\8,-1)\cup(0,+\8)$, $mrp\beta<0$, $m(m-p)>0$,
$\beta m[m-(\beta+1)p]<0$ and $(t-r)X_1<1\leqslant(t-r)X_2$,
\item[(B)] $\beta\in(-1,0)$, $pm\geqslant p^2$ or $\beta p^2<pm\leqslant0$ and $0<(t-r)X_2<1$,
\item[(C)] $\beta\in(-1,0)$, $r<0$, $(\beta+1)p^2<pm\leqslant p^2$ and $0<(t-r)X_1<1\leqslant(t-r)X_2$,
\item[(D)] $\beta\in(-\8,-1)\cup(0,+\8)$, $m(m-p)\leqslant0$, $\beta m[m-(\beta+1)p]<0$
and $0<(t-r)X_2<1$,
\eenu
$f(x)$ can not stay non-positive or non-negative on the whole interval $[0,T]$.
As a result, we can not directly use Jensen's inequality,
but the Semiconcave-semiconvex Theorem brings us some hope.

However, as $n$ increases or the parameters $\beta,t,r,m,p$ take general values,
the difficulty also increases greatly.
Therefore, we only present some concrete examples for these cases as follows.

\bx Under the conditions $r<0$ and (B), we let $n=4$, $m=\beta=-1/2$, $p=2$, $t=1$ and $r=-3$.
Then \eqref{3.zz1} and \eqref{3.1} hold, i.e.,
\be\label{3.51}\sum_{i=1}^4\(\frac{s+3a_i^2}{a_i}\)^{\frac12}\geqslant
2\sqrt{14}s^{\frac14}\geqslant\sqrt{7}\sum_{i=1}^4a_i^{\frac12},\ee
where $s=a_1^2+a_2^2+a_3^2+a_4^2$.
\ex

\bo Since $\disp\frac{m-\beta p}{p}=0.5\in(0,1)$, we can use Jensen's inequality and
obtain the second inequality of \eqref{3.51}.
In the following, we only consider the first inequality of \eqref{3.51}.

According to Theorem \ref{th3.1}, we know that
$$(t-r)X_2=\frac{2(\sqrt{6}-1)}{3}\in\(0,1\).$$
And hence $g(x)$ is convex on $\disp\left(0,\frac{X_2 s}{1-3X_2}\right]$ and concave
on $\disp\left[\frac{X_2 s}{1-3X_2},s\right)$.
We take arbitrarily a positive $\ve<\disp\min\left\{\frac{X_2}{1-3X_2},1-\frac{X_2}{1-3X_2}\right\}$
and set $a_1\leqslant a_2\leqslant a_3\leqslant a_4$.
Denote the left hand side of \eqref{3.51} by $F(a_1,a_2,a_3,a_4)$
with $a_1,a_2,a_3,a_4\in[\ve\sqrt{s},(1-\ve)\sqrt{s}]$.
By Theorem \ref{th2.2}, we know that $F(a_1,a_2,a_3,a_4)$
achieves its possible minimum in four cases in the following.

The first case is that $a_1^2=a_2^2=a_3^2=x^2$ and $a_4^2=s-3x^2\geqslant x^2$
with $x\in[\ve\sqrt{s},\sqrt{s}/2]$.
Let $\xi=x^2/s\in[\ve^2,1/4]$ and
$$h(\xi):=\frac1{\sqrt[4]{s}}F(a_1,a_2,a_3,a_4)
=\frac{3\sqrt{s+3x^2}}{\sqrt[4]{sx^2}}+\frac{\sqrt{4s-9x^2}}{\sqrt[4]{s(s-3x^2)}}
=\frac{3(1+3\xi)^{\frac12}}{\xi^{\frac14}}+\frac{(4-9\xi)^{\frac12}}{(1-3\xi)^{\frac14}}.$$
Then
$$h'(\xi)=-\frac34(1+3\xi)^{-\frac12}(1-3\xi)^{-\frac14}\left[\(\frac{1-3\xi}{\xi}\)^{\frac54}
+\frac{2-9\xi}{1-3\xi}\(\frac{1+3\xi}{4-9\xi}\)^{\frac12}\right].$$
Let $\eta=(1-3\xi)/\xi\in[1,+\8)$ and
$$\tilde{h}(\eta)=\(\frac{1-3\xi}{\xi}\)^{\frac54}
+\frac{2-9\xi}{1-3\xi}\(\frac{1+3\xi}{4-9\xi}\)^{\frac12}
=\eta^{\frac54}+\frac{2\eta-3}{\eta}\(\frac{\eta+6}{4\eta+3}\)^{\frac12}.$$
It is easy to see that $\tilde{h}(\eta)>0$ when $\eta\in[3/2,+\8)$.
When $\eta\in[1,3/2]$, we see that
\begin{align*}&\eta^{\frac52}-\(\frac{2\eta-3}{\eta}\)^2\frac{\eta+6}{4\eta+3}\\
=&\frac{1}{\eta^2(4\eta+3)}\left[(\eta-1)\(7\eta^4+7\eta^3+3\eta^2+9(6-\eta)\)
+\eta^{\frac92}(\eta^{\frac12}-1)(4\eta^{\frac12}-3)\right]\geqslant0,\end{align*}
which implies that $\tilde{h}(\eta)\geqslant0$ for all $\eta\in[1,+\8)$
and $h'(\xi)\leqslant0$ for all $\xi\in[\ve^2,1/4]$.
Hence
$$h(\xi)\geqslant h\(\frac14\)=2\sqrt{14}\hs\Ra\eqref{3.51}.$$

The second case is $a_1^2=a_2^2$, $a_4^2=(1-\ve)^2s$.
Since
$$a_1^2<\frac12(s-a_4^2)=\ve s-\frac12\ve s^2<\ve s,$$
we know
$$\frac1{\sqrt[4]{s}}F(a_1,a_2,a_3,a_4)>\frac1{\sqrt[4]{s}}\(\frac{s+3a_1^2}{a_1}\)^{\frac12}
>\frac1{\sqrt[4]{s}}\(\frac{s}{a_1}\)^{\frac12}
>\ve^{-\frac14}>2\sqrt{14},$$
if $\ve$ is sufficiently small.
Thus, we fix a sufficiently small $\ve_0\in(0,1/5)$ such that
$$F(a_1,a_2,a_3,(1-\ve)\sqrt{s})>2\sqrt{14}s^{\frac14}.$$

The third case is $a_3=a_4=(1-\ve_0)\sqrt{s}$.
In this case $a_3^2+a_4^2=2(1-\ve_0)^2s>s$, which is impossible.
The fourth case is $a_2=a_3=a_4=(1-\ve_0)\sqrt{s}$, which can be excluded either.
Eventually, we have proved \eqref{3.51} now.
\eo

\bx Under the conditions $r<0$ and (C), we let $n=4$, $\beta=-2/3$, $m=2/3$, $p=t=1$ and $r=-1$.
Then \eqref{3.zz2} and \eqref{3.2} hold, i.e.,
\be\label{3.52}\sum_{i=1}^4\left[(s+a_i)a_i\right]^{\frac23}\leqslant
\sqrt[3]{\frac{25s^4}4}\leqslant\sqrt[3]{25}\sum_{i=1}^4a_i^{\frac43},\ee
where $s=a_1+a_2+a_3+a_4$.
\ex

\bo The second inequality of \eqref{3.52} is obviously correct.
In this case, we have
$$(t-r)X_1=2(2-\sqrt3)\in(0,1).$$
Hence $g(x)$ is concave on $\disp\(0,\frac{X_1s}{1-X_1}\right]$
and convex on $\disp\left[\frac{X_1s}{1-X_1},s\)$.
Set $a_1\leqslant a_2\leqslant a_3\leqslant a_4$.
Denote the left hand side of \eqref{3.52} by $F(a_1,a_2,a_3,a_4)$
with $a_1,a_2,a_3,a_4\in[0,s]$.
By Theorem \ref{th2.2}, we know that $F(a_1,a_2,a_3,a_4)$
achieves its possible maximum in four cases in the following.

The first case is $a_1=a_2=a_3=x$ and $a_4=s-3x\geqslant x$ with $x\in[0,s/4]$.
Let $\xi=x/s\in[0,1/4]$ and
$$h(\xi)=s^{-\frac43}F(a_1,a_2,a_3,a_4)=3[(1+\xi)\xi]^{\frac23}+[(2-3\xi)(1-3\xi)]^{\frac23}.$$
Then
$$h'(\xi)=2(1+\xi)^{-\frac13}(1-3\xi)^{-\frac13}(2\xi+1)\left[\(\frac{1-3\xi}{\xi}\)^{\frac13}
+\frac{3(2\xi-1)}{2\xi+1}\(\frac{1+\xi}{2-3\xi}\)^{\frac13}\right]$$
Still we let $\eta=(1-3\xi)/\xi\in[1,+\8)$ and
$$\tilde{h}(\eta)=\(\frac{1-3\xi}{\xi}\)^{\frac13}
+\frac{3(2\xi-1)}{2\xi+1}\(\frac{1+\xi}{2-3\xi}\)^{\frac13}
=\eta^{\frac13}-\frac{3(\eta+1)}{\eta+5}\(\frac{\eta+4}{2\eta+3}\)^{\frac13}.$$
Since
$$\eta-\frac{27(\eta+4)}{2\eta+3}\(\frac{\eta+1}{\eta+5}\)^3
=\frac{2(\eta-1)(\eta^4+4\eta^3+7\eta^2+42\eta+54)}{(2\eta+3)(\eta+5)^3}\geqslant0,$$
we can conclude that $h'(\xi)>0$ on $(0,1/4]$ and $h(\xi)$ is strictly increasing
on $[0,1/4]$, which implies that
$$h(\xi)\leqslant h\(\frac14\)=\sqrt[3]{\frac{25}{4}}\hs\Ra\hs\eqref{3.52}.$$

The second case is $a_1=a_2$ and $a_4=s$, in which case $a_1=a_2=a_3=0$ and
$$F(a_1,a_2,a_3,a_4)=F(0,0,0,s)=2^{\frac23}s^{\frac43}<\sqrt[3]{\frac{25s^4}4}.$$
The third case is $a_3=a_4=s$, and the fourth case is $a_2=a_3=a_4=s$.
Both of the two cases are impossible.
As a result, we conclude \eqref{3.52}.
\eo

\bx Under the conditions $t>r>0$ and (D), we let $n=3$, $m=\beta\in(0,1)$, $p=r=1$ and $t=2$.
Then there is $\beta_0\in(0.5,1)$ such that when $\beta\in(0,\beta_0)$, \eqref{3.zz2} holds, i.e.,
\be\label{3.53}\(\frac{a_1}{a_1+2a_2+2a_3}\)^\beta+\(\frac{a_2}{2a_1+a_2+2a_3}\)^\beta
+\(\frac{a_3}{2a_1+2a_2+a_3}\)^\beta\leqslant\frac{3}{5^\beta}.\ee
\ex
\bo Following the proof of Theorem \ref{th3.1}, we know that $X_2=(1-\beta)/(1+\beta)$ and
hence $g(x)$ is concave on $(0,(1-\beta)s]$ and convex on $[(1-\beta)s,s)$.
We can as well set $a_1\leqslant a_2\leqslant a_3$.
Then by Theorem \ref{th2.2}, we know the left hand side of \eqref{3.53},
denoted by $F(a_1,a_2,a_3)$ for writing convenience,
achieves its possible maximum in three cases as follows.

The first case is $a_1=a_2=x$ and $a_3=s-2x\geqslant x$.
Then $x\in[0,s/3]$ and
$$F(a_1,a_2,a_3)=2\(\frac{x}{2s-x}\)^{\beta}+\(\frac{s-2x}{s+2x}\)^{\beta}.$$
We let $\xi=x/(2s-x)\in[0,1/5]$ and
$$h(\xi)=F(a_1,a_2,a_3)=2\xi^\beta+\(\frac{1-3\xi}{1+5\xi}\)^{\beta},$$
$$h'(\xi)=\frac{2\beta}{\xi^{1-\beta}(1+5\xi)^{1+\beta}}
\left[(1+5\xi)^{1+\beta}-4\(\frac{\xi}{1-3\xi}\)^{1-\beta}\right].$$
Now let $\eta=1+5\xi\in[1,2]$ and we can see that
$$(1+5\xi)^{\frac{1+\beta}{1-\beta}}-\frac{4^{\frac1{1-\beta}}\xi}{1-3\xi}
=\frac{3\eta^{\frac2{1-\beta}}-8\eta^{\frac{1+\beta}{1-\beta}}
+4^{\frac1{1-\beta}}\eta-4^{\frac1{1-\beta}}}{3\eta-8}.$$
Setting
$$\tilde{h}(\eta)=3\eta^{\frac2{1-\beta}}-8\eta^{\frac{1+\beta}{1-\beta}}
+4^{\frac1{1-\beta}}\eta-4^{\frac1{1-\beta}},$$
we have $\tilde{h}(1)=-5$, $\tilde{h}(2)=0$,
$$\tilde{h}'(\eta)=\frac{6}{1-\beta}\eta^{\frac{1+\beta}{1-\beta}}
-\frac{8(1+\beta)}{1-\beta}\eta^{\frac{2\beta}{1-\beta}}+4^{\frac{1}{1-\beta}},$$
\be\label{3.54}\tilde{h}''(\eta)=\frac{6(1+\beta)}{(1-\beta)^2}\eta^{\frac{2\beta}{1-\beta}}
-\frac{16\beta(1+\beta)}{(1-\beta)^2}\eta^{\frac{3\beta-1}{1-\beta}}
=\frac{2(1+\beta)}{(1-\beta)^2}\eta^{\frac{3\beta-1}{1-\beta}}\(3\eta-8\beta\),\ee
$$\tilde{h}'(1)=4^{\frac{1}{1-\beta}}-\frac{2(1+4\beta)}{1-\beta}
\hs\mb{and}\hs\tilde{h}'(2)=\frac{2-3\beta}{1-\beta}\cdot2^{\frac{2}{1-\beta}}.$$
We can see from \eqref{3.54} that only when $\beta\in(3/8,3/4)$, $h'(\eta)$ can
reach its least value in $(1,2)$, i.e.,
\be\label{3.55}\min_{\eta\in[1,2]}\tilde{h}'(\eta)=\tilde{h}'\(\frac{8\beta}{3}\)
=2^{\frac{2}{1-\beta}}-\frac{9}{8\beta^2}\(\frac{8\beta}{3}\)^{\frac{2}{1-\beta}}.\ee

Actually, in this process, we have to require $\tilde{h}'(\eta)\geqslant0$ on $[1,2]$,
which implies $\tilde{h}$ is non-decreasing on $[1,2]$ and so is $h'(\xi)$.
With these result, it yields that
$$\max_{\xi\in[1,1/5]}h(\xi)=h\(\frac15\)=\frac{3}{5^{\beta}}.$$
To this end, we only need to require
\be\label{3.56}\tilde{h}'(1)\geqslant0\mb{ when }\beta\in(0,3/8],\hs\tilde{h}'(2)\geqslant0
\mb{ when }\beta\in[3/4,1)\ee
\be\label{3.57}\mb{and}\hs\tilde{h}'\(\frac{8\beta}{3}\)\geqslant0,\mb{ when }\beta\in(3/8,3/4).\ee
It is not hard to deduce from \eqref{3.56} that $\beta\in(0,3/8]$.
From \eqref{3.55} and \eqref{3.57}, we get
$$\frac{8\beta^2}{9}\geqslant\(\frac{4\beta}{3}\)^{\frac{2}{1-\beta}}\hs\Lra\hs
\frac12\geqslant\(\frac{4\beta}{3}\)^{\frac{2\beta}{1-\beta}}\hs\Lra\hs
2\log_2\frac{3}{\beta}\geqslant\frac1\beta+3.$$
Let $\al=1/\beta\in(4/3,8/3)$ and
$$j(\al)=2\log_2(3\al)-\al-3\hs\mb{and so}\hs j'(\al)=\frac{2}{\al\ln 2}-1>0.$$
Since $j(4/3)=-1/3$ and $j(8/3)=1/3$, we can see that there is a unique $\al_0\in(4/3,8/3)$
such that $j(\al_0)=0$ and $j(\al)\geqslant0$ when $\al\in(\al_0,8/3)$.
By calculation using computers, we find
$$\beta_0=1/\al_0=0.5887287\cdots\in(0.5,1).$$
And hence it follows from \eqref{3.57} that $\beta\in(3/8,\beta_0)$.
As a result, we deduce that in this case when $\beta\in(0,\beta_0]$,
\eqref{3.53} holds true.

The second case is $a_3=s$. Then $a_1=a_2=0$ and $F(0,0,s)=1<3/5^{\beta_0}$.
The third case is $a_2=a_3=s$, which is impossible.
Consequently, we have obtained \eqref{3.53}.
\eo
In the example above, it has been proved that although the parameters satisfy
the condition (C), the conclusions of Theorems \ref{th3.1}, \ref{th3.2} and \ref{th3.3}
need not always hold.
Actually, under the cases (A), (B) and (C), it is still possible that none of the inequalities
\eqref{3.zz1}, \eqref{3.zz2}, \eqref{3.1} and \eqref{3.2} is valid.
The following example gives us a counterexample.

\bx Under the conditions $r<0$ and (A), we let $n=3$, $m=\beta=3/2$, $p=t=1$, $r=-1$ and
$$F(a_1,a_2,a_3)=\(\frac{a_1}{2a_1+a_2+a_3}\)^{\frac32}+\(\frac{a_2}{a_1+2a_2+a_3}\)^{\frac32}
+\(\frac{a_3}{a_1+a_2+2a_3}\)^{\frac32},$$
in which case $(t-r)X_1=2/5<1<2=(t-r)X_2$.
Then the right hand sides of \eqref{3.zz1}, \eqref{3.zz2}, \eqref{3.1} and \eqref{3.2} are
the same, i.e., $3/8$.
However, we observe that
$$\lim_{x\ra0^+}F(x,x,1)=F(0,0,1)=\frac{1}{2\sqrt2}<\frac38\hs\mb{and}\hs
\lim_{x\ra0^+}F(x,1,1)=F(0,1,1)=\frac{2}{3\sqrt3}>\frac38.$$
This situation thereby conflicts with all the inequalities \eqref{3.zz1},
\eqref{3.zz2}, \eqref{3.1} and \eqref{3.2}.
\ex

\section{Applications on inequality questions}

We employ the theorems in Section \ref{s3} to prove some interesting examples
and some mathematical competition questions in this section.
\subsection{Extensions on some inequalities}
The first example is a dimensional generalization of Example 7.19 of \cite{H11}.
This consequence also includes the result of Corollary 2.2 of \cite{WQB13}.
\bx Suppose that $a_1$, $\cdots$ $a_n>0$, $n\in\N^+$ and $n\geqslant2$.
Let
$$S_\beta=S_\beta(a_1,\cdots,a_n):=\sum_{i=1}^n\(\frac{a_i}{s-a_i}\)^\beta$$
\be\label{4.0}\mb{and}\hs \beta_n=\frac{\ln n-\ln(n-1)}{\ln(n-1)-\ln(n-2)}\hs\mb{for }n>2.\ee
where $s=a_1+\cdots+a_n$. Then $\beta_n$ is increasing in $n$ and
\be\label{4.1}\inf S_\beta(a_1,\cdots,a_n)=
\left\{\ba{cl}2,&\beta\in(0,\beta_3);\\
\disp\frac{k}{(k-1)^{\beta}},&\beta\in[\beta_k,\beta_{k+1}),\;k=3,\cdots,n-1;\\[2ex]
\disp\frac{n}{(n-1)^\beta},&\beta\in(-\8,0]\cup[\beta_n,+\8).\ea\right.\ee
\ex
\bo First by definition \eqref{4.0}, using Cauchy mean value theorem, we know that for each $n>2$,
there is $\theta_n\in(0,1)$ such that
$$\beta_n=\frac{1/(n-1+\theta_n)}{1/(n-2+\theta_n)}=\frac{n-2+\theta_n}{n-1+\theta_n}.$$
Hence $\beta_n$ is increasing as $n$ increases.
Moreover,
$$\lim_{n\ra2^+}\beta_n=0\hs\mb{and}\hs \lim_{n\ra+\8}\beta_n=1,$$
and
$$\frac{n-1}{(n-2)^{\beta_n}}=\frac{n}{(n-1)^{\beta_n}}\hs\mb{and}\hs\beta_n>\frac{n-2}{n}.$$
And hence
\be\label{5.53C}\mb{when }\beta\in[\beta_n,1),\hs
\frac{n-1}{(n-2)^{\beta}}\geqslant\frac{n}{(n-1)^{\beta}};\ee
\be\label{5.53D}\mb{when }\beta\in\(0,\beta_n\),\hs
\frac{n-1}{(n-2)^{\beta}}<\frac{n}{(n-1)^{\beta}}.\ee
\Vs

Next we show \eqref{4.1}.
In accord with Theorem \ref{th3.1},
the conditions presented in this example is $p=t=r=1$ and $m=\beta$.
Then by the cases (ii.1) and (iii.3) of Theorem \ref{th3.1}
and (ii.2) of Theorem \ref{th3.2},
we know that if $\beta\geqslant1$ or $\beta\leqslant0$, \eqref{4.1} can be deduced
(it is obvious when $\beta=0$).
Next we only consider the case when $\beta\in(0,1)$.

Following the proof of Theorem \ref{th3.1}, we know that
$$g''(x)=\frac{\beta sx^{\beta-2}}{(s-x)^{\beta+2}}\left[2x-(1-\beta)s\right],
\hs x\in[0,s].$$
This means that $g(x)=(x/(s-x))^\beta$ is concave on $[0,(1-\beta)s/2]$ and
convex on $[(1-\beta)s/2,s]$.
Let $a_k$ be non-decreasing when $k$ increases.
By Theorem \ref{th2.2}, we pick one arbitrary possible minimum point of $S_\beta(a_1,\cdots,a_n)$
(Here we allow $a_k=0$ for $k\leqslant n-2$) such that
$$a_1=\cdots=a_{n-k}=0<a_{n-k+1}=s-(k-1)x\leqslant a_{n-k+2}=\cdots=a_n=x,$$
with $k\geqslant2$.
It is easy to see that $x\in[s/k,s/(k-1))$.
We let $\xi=x/s\in[1/k,1/(k-1))$ and
$$h(\xi)=S_\beta(a_1,\cdots,a_n)=\(\frac{1-(k-1)\xi}{(k-1)\xi}\)^{\beta}
+(k-1)\(\frac{\xi}{1-\xi}\)^{\beta}.$$
We know that
\be\label{3.53X}\mb{if }k=2,\hs
h(\xi)\geqslant2\hs\mb{and}\hs h(\xi)=2\mb{ if and only if }\xi=\frac12.\ee

In the following, we consider the case when $k\geqslant3$.
Let $\disp\eta=\frac{1-(k-1)\xi}{(k-1)\xi}\in\(0,\frac1{k-1}\right]$ and
$$i(\eta):=h(\xi)=\eta^\beta+\frac{k-1}{[(k-1)\eta+k-2]^{\beta}}.$$
Then
$$i'(\eta)=\frac{\beta}{[(k-1)\eta+k-2]^{\beta+1}}
\left[[(k-1)\eta+k-2]^{\beta+1}\eta^{\beta-1}-(k-1)^2\right].$$
Let $j(\eta):=[(k-1)\eta+k-2]^{\beta+1}\eta^{\beta-1}-(k-1)^2$ and then
$$j'(\eta)=[(k-1)\eta+k-2]^{\beta}\eta^{\beta-2}[2\beta(k-1)\eta+(\beta-1)(k-2)].$$
Setting $\disp\eta_0=\frac{(1-\beta)(k-2)}{2\beta(k-1)}$,
we split it into two situations for discussing.

If $\eta_0\geqslant\disp\frac1{k-1}$, i.e., $\disp\beta\in\(0,\frac{k-2}{k}\right]$, then
$$j'(\eta)\leqslant0,\hs j(\eta)\geqslant j\(\frac1{k-1}\)=0,\hs i'(\eta)\geqslant0$$
\be\label{4.2}\mb{and}\hs h(\xi)=i(\eta)>\lim_{\eta\ra0^+}i(\eta)=\frac{k-1}{(k-2)^\beta}.\ee

If $\disp\eta_0\in\(0,\frac1{k-1}\)$, i.e., $\disp\beta\in\(\frac{k-2}{k},1\)$, then
$$j'(\eta)<0\mb{ as }\eta\in(0,\eta_0),\hs j'(\eta_0)=0,\hs j'(\eta)>0
\mb{ as }\eta\in\(\eta_0,\frac1{k-1}\).$$
Then $j$ reaches its minimum at $\eta=\eta_0$.
Noting that
$$\lim_{\eta\ra0^+}j(\eta)=+\8\hs\mb{and}\hs j\(\frac1{k-1}\)=0,$$
we know there exists $\eta_1\in(0,\eta_0)$ such that
$$j(\eta)>0\mb{ as }\eta\in(0,\eta_1),\hs j(\eta_1)=0,\hs j(\eta)<0
\mb{ as }\eta\in\(\eta_1,\frac1{k-1}\).$$
This also means that
$$i'(\eta)>0\mb{ as }\eta\in(0,\eta_1),\hs i'(\eta_1)=i'\(\frac1{k-1}\)=0,\hs i'(\eta)<0
\mb{ as }\eta\in\(\eta_1,\frac1{k-1}\).$$
Hence $i$ is increasing on $(0,\eta_1]$ and decreasing on $\disp\left[\eta_1,\frac1{k-1}\right]$.
Therefore, recalling \eqref{4.2}, we obtain for all $\beta\in(0,1)$,
\be\label{3.53B}S_\beta(a_1,\cdots,a_n)\geqslant\min\left\{\lim_{\eta\ra0^+}i(\eta),
i\(\frac1{k-1}\)\right\}
=\min\left\{\frac{k-1}{(k-2)^{\beta}},\frac{k}{(k-1)^{\beta}}\right\}.\ee

Now we consider $n\geqslant 2$ and combine \eqref{3.53X}, \eqref{3.53B}, \eqref{5.53C}
and \eqref{5.53D} to obtain the following consequences.
When $\beta\in[\beta_n,1)$, by \eqref{3.53X}, \eqref{3.53B} and \eqref{5.53C}, we know
\be\label{4.7}\inf S_\beta(a_1,\cdots,a_n)=\min\left\{2,\frac{3}{2^\beta},\frac{4}{3^{\beta}},\cdots,
\frac{n}{(n-1)^{\beta}}\right\}=\frac{n}{(n-1)^{\beta}}.\ee
When $\beta\in[\beta_k,\beta_{k+1})$, $k=3,\cdots,n-1$, by \eqref{3.53X},
\eqref{3.53B}, \eqref{5.53C} and \eqref{5.53D}, we know
\be\label{4.8}\inf S_\beta(a_1,\cdots,a_n)=\min\left\{2,\frac{3}{2^\beta},\frac{4}{3^{\beta}},\cdots,
\frac{n}{(n-1)^{\beta}}\right\}=\frac{k}{(k-1)^{\beta}}.\ee
When $\beta\in(0,\beta_3)$, by \eqref{3.53X},
\eqref{3.53B}, \eqref{5.53C} and \eqref{5.53D}, we know
\be\label{4.9}\inf S_\beta(a_1,\cdots,a_n)=\min\left\{2,\frac{3}{2^\beta},\frac{4}{3^{\beta}},\cdots,
\frac{n}{(n-1)^{\beta}}\right\}=2.\ee
At last, \eqref{4.1} follows from \eqref{4.7}, \eqref{4.8} and \eqref{4.9}.
The proof is complete.
\eo

Based on Theorem \ref{th3.1} in Section \ref{s3}, we can also get a more general result
as follows.
The following example is a new generalization of Mitrinovi\'c inequality (see \cite{FY22}).
\bx\label{x4.2} Under the conditions of Theorem \ref{th3.1}, we further pick $k\in\N^+\cap[1,n-1]$.
Then in the cases (i), (ii), (iii) and (iv) with $t>r$ replaced by $t>kr$ where $t>r$ appears,
\be\label{4.11}\sum_{i=1}^n\frac{(a^p_i+\cdots+a^p_{i+k-1})^m}{[ts-r(a^p_i+\cdots+a^p_{i+k-1})]^\beta}
\geqslant
\frac{k^mn^{\beta+1-m}}{(nt-kr)^{\beta}}\(\sum_{i=1}^na_i^p\)^{m-\beta};\ee
in the cases (v), (vi), (vii) and (viii) with $t>r$ replaced by $t>kr$ where $t>r$ appears,
\be\label{4.12}\sum_{i=1}^n\frac{(a^p_i+\cdots+a^p_{i+k-1})^m}{[ts-r(a^p_i+\cdots+a^p_{i+k-1})]^\beta}
\leqslant
\frac{k^mn^{\beta+1-m}}{(nt-kr)^{\beta}}\(\sum_{i=1}^na_i^p\)^{m-\beta}.\ee
Here $a_{k}$ is supposed to be $a_{k-n}$ if $k>n$.
\ex
\bo We only show \eqref{4.11} since \eqref{4.12} can be proved similarly.
Let $A_i:=a_i^p+\cdots+a_{i+k-1}^p$ for $i=1,\cdots,n$ and $S:=A_1+\cdots+A_n=ks$.
Since $t>kr$ implies $t/k>r$, then under the cases (i), (ii), (iii) and (iv),
we can use Theorem \ref{th3.1} and obtain
\begin{align*}&\sum_{i=1}^n\frac{(a^p_i+\cdots+a^p_{i+k-1})^m}{[ts-r(a^p_i+\cdots+a^p_{i+k-1})]^\beta}\\
=&\sum_{i=1}^n\frac{A_i^m}{\disp\(\frac{t}{k}S-rA_i\)^\beta}
\geqslant\frac{n^{\beta+1-m}}{\disp\(\frac{nt}{k}-r\)^{\beta}}S^{m-\beta}
=\frac{k^mn^{\beta+1-m}}{(nt-kr)^{\beta}}\(\sum_{i=1}^na_i^p\)^{m-\beta},\end{align*}
which ends the proof.\eo

The theorems in Section \ref{s3} can be used to prove inequalities concerning with
the sides of triangles.
\bx Let $a$, $b$ and $c$ be three sides of a triangle. Then
\be\label{4.13}\frac{a^2}{b+c-a}+\frac{b^2}{c+a-b}+\frac{c^2}{a+b-c}\geqslant a+b+c.\ee
\ex
\bo Since $a$, $b$ and $c$ are the sides of a triangle,
the sum of each two of $a$, $b$ and $c$ is bigger than the other.
Therefore, we can take $n=3$, $m=2$, $p=t=\beta=1$ and $r=2$ in \eqref{3.zz1}
in the case (iii.3) of Theorem \ref{th3.1} and directly obtain \eqref{4.13}.
\eo

\subsection{Applications on competition questions}
In the following, we present some mathematical competition questions that can be obtained by
the main theorems in Section \ref{s2}.
For writing convenience, we denote the left hand side of some inequality by LHS.

\bx[28th IMO Pre-selection Question]
Let $a$, $b$ and $c$ be the sides of a triangle and $2S=a+b+c$.
Prove that
$$\frac{a^m}{b+c}+\frac{b^m}{c+a}+\frac{c^m}{a+b}\geqslant\(\frac23\)^mS^{m-1},$$
where $m\geqslant1$.
Particularly, when $m=2$, this is a question of 19th Nordic Mathematical Olympiad Contest in 2005.
\ex
\bo It is a simple example of \eqref{3.zz1}
by picking $n=3$, $m\geqslant2$ and $p=\beta=t=r=1$ in the case (iii.3) of Theorem \ref{th3.1},
and an example of \eqref{3.1} by picking $n=3$ and $m=p=\beta=t=r=1$ in the case (i.1) of Theorem \ref{th3.2}
(This is also the famous Nesbitt's inequality).
\eo

\bx[31st IMO Pre-selection Question] Let $a$, $b$, $c$ and $d$
be positive real numbers such that $ab+bc+cd+da=1$.
Prove
$$\frac{a^3}{b+c+d}+\frac{b^3}{c+d+a}+\frac{c^3}{d+a+b}+\frac{d^3}{a+b+c}\geqslant\frac13.$$
This example was also selected in Chinese Mathematical Olympiad in Senior
(Xinjiang Division) Preliminary Contest in 2020.
\ex
\bo Pick $n=4$, $m=3$ and $p=\beta=t=r=1$ in \eqref{3.1} in the case (i.1) of Theorem \ref{th3.2}.
Then we have
\begin{align*}&\mb{LHS}\geqslant\frac13(a^2+b^2+c^2+d^2)\\
=&\frac13\frac{(a^2+b^2)+(b^2+c^2)+(c^2+d^2)+(d^2+a^2)}{2}\\
\geqslant&\frac13(ab+bc+cd+da)=\frac13.
\end{align*}
The proof is hence over.\eo

\bx[IMO-36 in 1995]Let $a$, $b$, $c$ be positive real numbers such that $abc=1$.
Prove that
$$\frac1{a^3(b+c)}+\frac1{b^3(c+a)}+\frac1{c^3(a+b)}\geqslant\frac32.$$
\ex
\bo Picking $n=3$, $m=-2$, $p=-1$ and $\beta=t=r=1$ in \eqref{3.zz1} in the case (iii.3)
of Theorem \ref{th3.1}, we obtain
\begin{align*}&\mb{LHS}=\frac{bc}{a^2(b+c)}+\frac{ca}{b^2(c+a)}+\frac{ab}{c^2(a+b)}\\
=&\frac1{a^2(b^{-1}+c^{-1})}+\frac1{b^2(c^{-1}+a^{-1})}+\frac1{c^2(a^{-1}+b^{-1})}
\geqslant\frac12(a^{-1}+b^{-1}+c^{-1})\\
\geqslant&\frac32\sqrt[3]{(abc)^{-1}}=\frac32,\end{align*}
where we have also used the mean value inequality.
\eo

\bx[Serbian Math Olympiad in 2005]Let $x$, $y$ and $z$ be positive numbers.
Prove
$$\frac{x}{\sqrt{y+z}}+\frac{y}{\sqrt{z+x}}+\frac{z}{\sqrt{x+y}}
\geqslant\sqrt{\frac32(x+y+z)}.$$
\ex
\bo Pick $n=3$, $m=p=t=r=1$ and $\beta=1/2$ in \eqref{3.zz1} in the case (iii.3)
of Theorem \ref{th3.1}, we easily obtain the conclusion.
\eo

\section{Applications on Hurwitz-Lerch zeta functions}
The Hurwitz-Lerch zeta function $\zeta(z,\beta,a)$ is defined by
$$\zeta(z,\beta,a):=\sum_{n=0}^{\8}\frac{z^n}{(n+a)^\beta},$$
where $a\in\C\sm\Z_0^-$, $\beta\in\C$ when $|z|<1$ and $\fR(\beta)>1$ when $|z|=1$.
Here $\C$ is the set of complex numbers, $\Z_0^-$ is the set of nonpositive integers
and $\fR(\beta)$ means the real part of $\beta\in\C$.

In the following theorem, we only discuss the relation about Hurwitz-Lerch zeta functions
with real variables.

\bt\label{th5.1} Let $x=(x_1,\cdots,x_n)\in\R^n$, $\{a_n\}_{n\in\N^+}$ be a positive sequence,
$$X=\{x\in\R^n:x_i>0\mb{ and }x_1+\cdots+x_n=1\},\hs\hs
\al=\inf_{n\in\N^+}a_n,$$
$z>0$ and $r\in\R$ such that $\al>r$.
Then in the case when $\beta r>0$,
\be\label{5.1}\min_{x\in X}\sum_{i=1}^{n}x_i\zeta\(z,\beta,a_n-rx_i\)
=\zeta\(z,\beta,a_n-\frac rn\);\ee
in the following cases
\benu\item[(1)] $\beta>0$, $r<0$ and $2(\al-r)\geqslant-(\beta+1)r$,
\item[(2)] $\beta<-1$, $r>0$ and $2(\al-r)\geqslant-(\beta+1)r$,
\item[(3)] $\beta\in[-1,0)$ and $r>0$,\eenu
\be\label{5.2}\max_{x\in X}\sum_{i=1}^{n}x_i\zeta\(z,\beta,a_n-rx_i\)
=\zeta\(z,\beta,a_n-\frac rn\).\ee
Here each Hurwitz-Lerch zeta function in \eqref{5.1} and \eqref{5.2}
is supposed to be convergent.
\et

\br In Theorem \ref{th5.1}, we only consider the case when $\beta r\ne0$,
since when $\beta r=0$, $\zeta(z,\beta,a_n-rx_i)$ does not depend on $x_i$,
and it is obvious that
$$\sum_{i=1}^{n}x_i\zeta\(z,\beta,a_n-rx_i\)=\zeta\(z,\beta,a_n-\frac rn\).$$
\er

\noindent\textit{Proof of Theorem \ref{th5.1}.} In this proof, we always assume that $m=p=1$.
For the case when $\beta r>0$, we first consider the case when $\beta>0$, $r>0$
or $\beta<-1$, $r<0$.
In this case, (iii.3) of Theorem \ref{th3.1} is satisfied.
Thus we obtain by \eqref{3.zz1} that
\be\label{5.3}\sum_{i=1}^n\frac{x_i}{(j+a_n-rx_i)^\beta}\geqslant
\frac{n^{\beta}}{[n(j+a_n)-r]^\beta}=\frac{1}{(j+a_n-\frac rn)^\beta}.\ee
Then multiplying \eqref{5.3} by $z^j$ and
adding the results for $j=0,1,\cdots,k$ together with $k\in\N^+$,
we have
\be\label{5.4}\sum_{i=1}^nx_i\sum_{j=0}^k\frac{z^j}{(j+a_n-rx_i)^\beta}=
\sum_{j=0}^kz^j\sum_{i=1}^n\frac{x_i}{(j+a_n-rx_i)^\beta}\geqslant
\sum_{j=0}^k\frac{z^j}{(j+a_n-\frac{r}n)^\beta}.\ee
Let $k$ tend to the infinity, we conclude that
\be\label{5.5}\sum_{i=1}^nx_i\zeta(z,\beta,a_n-rx_i)\geqslant\zeta(z,\beta,a_n-\frac rn),\ee
which implies \eqref{5.1}.
For the case when $\beta=-1$ and $r<0$, (ii.1)of Theorem \ref{th3.1} is satisfied;
For the case when $\beta\in(-1,0)$ and $r<0$, we see that
$$(j+a_n-r)X_2=\frac{2(j+a_n-r)}{(\beta+1)(-r)}>\frac{2}{\beta+1}>1$$
and (iv) is satisfied.
In these two cases, we also have \eqref{5.3} and hence \eqref{5.1}.

Next we show \eqref{5.2} for the case when $\beta r<0$.
When $\beta>0$, $r<0$ and $2(\al-r)\geqslant-(\beta+1)r$
(the proof for the case (2) is the same), we have
$$(j+a_n-r)X_2=-\frac{2(j+a_n-r)}{(\beta+1)r}\geqslant
-\frac{2(\al-r)}{(\beta+1)r}\geqslant1.$$
Thus (viii) of Theorem \ref{th3.1} is satisfied.
Similar to the discussion of \eqref{5.3}, \eqref{5.4} and \eqref{5.5},
we can use \eqref{3.zz2} to obtain \eqref{5.2}.
When $\beta=-1$ and $r>0$, (vi.1) of Theorem \ref{th3.1} is satisfied.
When $\beta\in(-1,0)$ and $r>0$, (vii.2) of Theorem \ref{th3.1} is satisfied.
As a result, \eqref{5.2} can be similarly obtained.
The proof is hence finished now.
\qed

\br Theorem \ref{th5.1} is a generalization of Theorem 3.2 of \cite{WQB13}.
Specifically, when $a_n=1+\frac1n$, $r=1$ and $z=1$, \eqref{5.1} is the result obtained
in \cite{WQB13}.
\er

\br This article mainly generalizes Nesbitt's inequality in respect of dimensions and parameters
and gives different results in various cases.
The argument also provides a series of methods to estimate algebraic expressions analogous
to \eqref{1.1}.
This article is not concerning with the inequalities with weights like
\cite{BS12,BS13-1,BS13-2,BZ08,BP11,WQB13}.
Actually, it is still interesting to study the inequalities \eqref{3.zz1}, \eqref{3.zz2},
\eqref{3.1} and \eqref{3.2} with weights.
\er

\section*{Acknowledgements}

Our work was supported by grant from the National Natural Science Foundation of China
(NSFC No. 11801190).


\end{document}